\documentclass[12pt]{amsart}

\usepackage{epsfig}
\usepackage{verbatim}
\usepackage{undertilde}

\begin{document}
\title {Unlabeled equivalence for matroids representable over finite fields}

\maketitle
\begin {center}
S. R. Kingan \\
Department of Mathematics \\
Brooklyn College, City University of New York\\
 Brooklyn, NY 11210\\
skingan@brooklyn.cuny.edu\\  
\end {center}
\bigskip

\begin{abstract}  We present a new type of equivalence for representable matroids that uses the automorphisms of the underlying matroid. Two $r\times n$ matrices $A$ and $A'$ representing the same matroid $M$ over a field  $F$ are  {\it geometrically equivalent representations} of $M$ if one can be obtained from the other by elementary row operations, column scaling, and column permutations. Using geometric equivalence, we give a method for exhaustively generating non-isomorphic matroids representable over a finite field $GF(q)$, where $q$ is a power of a prime.  
\end{abstract}

\section{Introduction}

This paper introduces a new type of equivalence for representable matroids that takes advantage of the automorphisms of the underlying matrix.  It is useful because with this type of equivalence it becomes no harder to handle a matroid  representable matroid over $GF(q)$, where $q$ is a prime power, than a matroid representable over $GF(2)$. In particular, for a fixed rank $r$, there is a polynomial time algorithm in $r$ with leading coefficient $q^r$ that determines whether or not two matrices over $GF(q)$ are geometrically equivalent.

The question of representability over a finite field is complicated by the presence of inequivalent representations. For example, the three matrices $A$, $B$, $C$ below, represent the same matroid $W^3$ (the 3-whirl). 

\tiny
\[ 
A=\left[ 
\begin{array}{ccc|ccc}
&&&   1&0&1 \\
&I_3&&1&1&0 \\
&&&   0&1&1
\end{array} 
\right] 
B=\left[ 
\begin{array}{ccc|ccc}
&&&   1&0&1 \\
&I_3&&1&1&0 \\
&&&   0&1&2
\end{array} 
\right] 
C=\left[ 
\begin{array}{ccc|ccc}
&&&   1&0&1 \\
&I_3&&1&1&0 \\
&&&   0&1&3
\end{array} 
\right] 
\] 
\normalsize

\noindent  Add the column $[1, 1, 1]^T$ to each and we see that $A$ with $[1, 1, 1]^T$ represents the matroid $F_7^-$ (non-Fano), whereas  $B$ and $C$ with $[1, 1, 1]^T$, respectively, represent the matroid $X_7$ in Figure 1. Moreover, there is no column that can be added to $A$ to get a matrix that represents $X_7$. So, if one happens to consider $A$ as a representation of $W^3$ over $GF(5)$, then there is no indication that $X_7$ is a single-element extension of $W^3$,  and the fact that $X_7$ is representable over $GF(5)$ could conceivably be missed.

The matroid terminology follows Oxley [4]. At present, there are two kinds of equivalence in representable matroids. Two $r\times n$ matrices $A$ and $A'$ representing the same matroid $M$ over a field  $F$ are said to be {\it projectively equivalent representations} of $M$ if one can be obtained from the other by elementary row operations or column scaling. Otherwise $A$ and $A'$ are {\it projectively inequivalent}. The matrices $A$ and $A'$ are said to be {\it (algebraically) equivalent} representations of $M$ if one can be obtained from the other by elementary row operations, column scaling, or field automorphisms.  Otherwise they are called {\it inequivalent representations}. A matroid is uniquely representable over $F$ if all of its representations are equivalent [3]. For prime fields, projective equivalence and algebraic equivalence coincide, as prime fields do not have non-trivial automorphisms. For fields of order $q$, where $q$ is a prime power, projective equivalence is a refinement of  equivalence. So, a result that holds for projective equivalence automatically holds for equivalence. 

We define 
$A$ and $A'$ to be  {\it geometric equivalent representations} of $M$ if one can be obtained from the other by elementary row operations, column scaling, or column permutations. Otherwise $A$ and $A'$ are {\it geometric inequivalent}.  This definition disregards the labels on the columns and may be viewed as an ``unlabeled" equivalence.  We call it geometric equivalence because allowing column permutations is the same allowing linear transformations. 

The three representations of $W^3$ over $GF(5)$, mentioned earlier, are projectively inequivalent representations. But there is a linear transformation that maps $B$ to $C$. Specifically, we can convert $B$ to $C$ by the following sequence of operations: swap rows 1 and 3; multiply column 6 by 3; and swap columns 1 and 3 ad columns 4 and 5. There is, however, no linear transformation that maps $A$ to $B$. So, $A$ and $B$ are geometrically inequivalent representations of $W^3$ over $GF(5)$.    

Besides doing an exhaustive search, which we did, the fact that there is no linear transformation that maps $A$ to $B$ can be confirmed by checking single-element extensions. Matrix $A$ has six non-isomorphic single-element extensions. Let us call them $A_7$, $B_7$, $C_7$, $D_7$, and the well-known $F_7^-$ and $P_7$. Matrix $B$  and $C$  each have seven non-isomorphic single-element extensions, $A_7$, $B_7$, $C_7$, $D_7$, $X_7$, $Y_7$ and $P_7$ (see Figure 1).

\begin{figure}[h]
\centering
\epsfxsize 5in \epsfbox{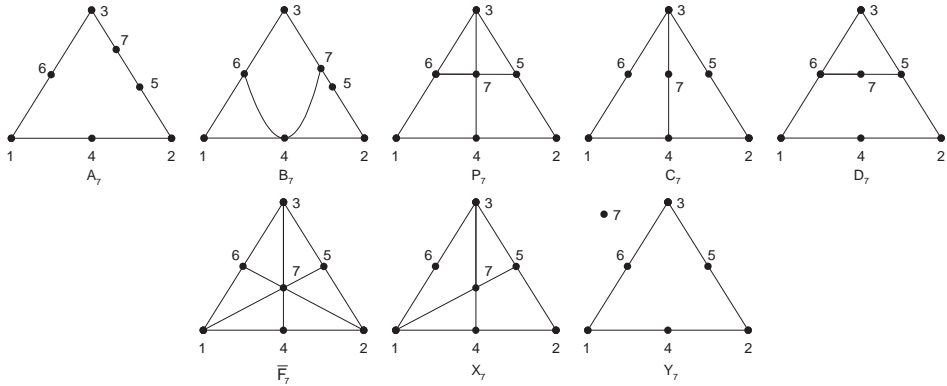}
\caption{Single-element extensions of $W^3$ over $GF(5)$ \label{ExtensionsOf3Whirl}}
\end{figure}
\bigskip 

In Section 2 we give a method for coordinatizing a matroid over a finite field taking into account geometric equivalence. In Section 3 we give a method for generating single-element extensions and coextensions using geometric equivalence. This gives an algorithm for isomorph-free and exhaustive generation of representable matroids.

\section{Method for coordinatizing a matroid over $GF(q)$}

In this section, we modify Brylawski and Lucas's method [1] for coordinatizing a matroid over a field [1] to take into account geometric equivalence. To obtain a matrix representation of a matroid $M$ with ground set $E$ and basis $B=\{1, \dots ,r\}$, first use the fundamental circuits of each element in $E-B$ with respect to $B$ to construct the matrix $[I_r |D^\sharp ]$, where column $k$ of $D^\sharp$ has ones corresponding to the fundamental circuit of $k$ with respect to $B$ and zeros elsewhere.

Observe that $D^\sharp$ can be viewed as an incidence matrix with the rows labeled by 
$B=\{1,  \dots , r\}$ and columns labeled by $E-B=\{r+1,  \dots , n\}$.  Let 
$G$ be the bipartite graph corresponding to this incidence matrix.  The vertices 
in the two classes of $G$ are labeled with elements of $B$ and 
$E-B$, respectively.  The edges in $G$ correspond to the ones in $D^\sharp$.  

Select any spanning forest $B_G$ of $G$ and assign the corresponding 
entries in $D^\sharp$ arbitrary values from $F$.  Since $B_G$ is a basis, we can do row and column scaling to reduce all 
these values to one.   So, without loss of generality, we may assume these entries are one. Circle them and label the remaining non-zero entries of $D^\sharp$ as $(a, b, c, \dots )$. Call the resulting matrix $D$.  The values of $(a, b, c, \dots )$ may be found by setting up a 
system of equations using the circuits of $M$ and solving the system over $F$.  If the 
system has a solution, then the matroid is representable over $F$. Otherwise it is {\it not representable} over $F$.  Two situations arise:

\begin {enumerate}
\item [(i)] If there is more than one value for the ordered sequence $(a, b, c, \dots  )$ over $F$, then the matrices $[I_r|D]$ with the different values of  $(a, b, c, \dots  )$ are projectively inequivalent representations  of $M$. If there is just one sequence $(a, b, c, \dots  )$ over $F$, then $M$ is uniquely representable over $F$ with respect to projective equivalence.
\item [(ii)] If  $(a_1, b_1, c_1, \dots  )$ and $(a_2, b_2, c_2, \dots  )$ are two ordered sequences over $F$ with corresponding matrices  $[I_r|D_1]$ and $[I_r|D_2]$ and there is a linear transformation that maps one matrix to the other, then the two matrices are geometrically equivalent. Otherwise they are   geometrically ineqivalent representations of $M$.
Geometric equivalence is an equivalence relation on the set of sequences $(a, b, c, \dots  )$ over $F$. Two sequences are geometrically equivalent if they belong to the same equivalence class. The number of geometrically inequivalent representations is the number of equivalence classes. If there is only one equivalence class, then $M$ is  uniquely representable over $F$ with respect to geometric equivalence.
\end {enumerate}
 
The next proposition follows from this discussion and  contrasts projective equivalence with geometric equivalence.  Without loss of generality, we state results only for connected matroids.

\bigskip
\noindent {\bf Proposition 2.1.}  {\it Let $A_1=[I_r|D_1]$ and $A_2=[I_r|D_2]$ be matrices representing the same connected matroid $M$ over $GF(q)$, such that columns of $A_1$ represent the same elements of $M$ as corresponding columns of $A_2$. Then $G_1=G_2$. Assume that $F$ is a spanning tree of $G$, having the property that, for each edge $f\in F$, the entries in $D_1$ and $D_2$ corresponding to $f$ are equal. Then 
\begin{enumerate}
\item[(i)] $A_1$ and $A_2$ are projectively equivalent if and only if $D_1 = D_2$.
\item[(ii)] $A_1$ and $A_2$ are geometrically equivalent if and only if there is a linear transformation from $A_1$ to $A_2$. 
\end{enumerate}
} 
\noindent Moreover, the number of geometrically inequivalent representations of $M$ is at most the number of projectively inequivalent representations of $M$.
$\qed$
\bigskip

We illustrate this method by coordinatizing the matroid $Q_6$ shown in Figure 2. The fundamental 
circuits corresponding to the basis $\{1, 2, 3\}$ are $\{1, 2, 4\}$, $\{2, 3, 5\}$, and $\{1, 2, 3, 6\}$.  The matrix obtained from the fundamental circuits is

\small 
\[ 
\left[ 
\begin{array}{ccc|ccc}
&&&1&0&1 \\
&I_3&&1&1&1 \\
&&&0&1&1
\end{array} 
\right] 
\] 
\normalsize

\begin{figure}[h]
\centering
\epsfxsize 2in \epsfbox{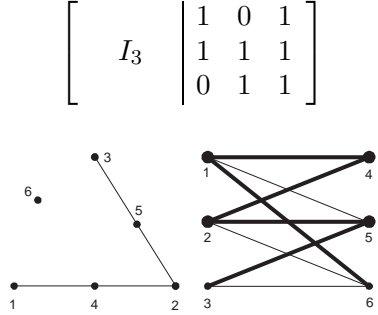}
\caption{Coordinatizing $Q_6$ \label{Q6}}
\end{figure}
\bigskip

\noindent The corresponding bipartite graph is shown in Figure 2 with a spanning 
tree highlighted.  The entries corresponding to the edges of the spanning tree 
may be taken as 1 and the remaining entries are denoted by $(a, b)$ where $a$ and $b$ are non-zero members of the field.  
Thus we can write the matrix representing $Q_6$  as

\small
\[ 
\left[ 
\begin{array}{ccc|ccc}
&&&1&0&1 \\
&I_3&&1&1&a \\
&&&0&1&b
\end{array} 
\right] 
\] 
\normalsize

\noindent  Over  $GF(3)$ the choices for the ordered sequence $(a, b)$ are $(1, 1)$, $(1, 2)$, $(2, 1)$ and $(2, 2)$. The first gives $W_3$ and the rest give $W^3$. So $Q_6$ is not representable over $GF(3)$. 
  Over $GF(5)$ the ordered sequence $(a, b)$ may be $(1, 1)$, $(1, 2)$, $(2, 1)$, $(1, 3)$, $(3, 1)$, $(1, 4)$, $(4, 1)$, $(2, 2)$, $(2, 3)$, $(3, 2)$, $(2, 4)$, $(4, 2)$, $(3, 3)$, $(3, 4)$, $(4, 3)$, or $(4, 4)$. Of these $(1, 1)$ gives $W_3$; $(1, 2)$, $(2, 1)$, $(1, 3)$, $(1, 4)$, $(2, 2)$, $(3, 2)$,  $(3, 3)$, $(4, 3)$, and $(4, 4)$ give $W^3$ and the remaining give $P_6$. Let us refer to the matrices with $(a, b)$ equal to $(3, 1)$, $(4, 1)$, $(2, 3)$, $(2, 4)$, $(4, 2)$ and $(3, 4)$ as $B_1$, $B_2$, $B_3$, $B_4$,  $B_5$, and $B_6$, respectively. So, $Q_6$ has six projectively inequivalent representations over $GF(5)$, shown below.

\tiny
\[ 
B_1=\left[ 
\begin{array}{c|ccc}
&1&0&1 \\
I_3&1&1&3 \\
&0&1&1
\end{array} 
\right] 
B_2=\left[ 
\begin{array}{c|ccc}
&1&0&1 \\
I_3&1&1&4 \\
&0&1&1
\end{array} 
\right] 
B_3=\left[ 
\begin{array}{c|ccc}
&1&0&1 \\
I_3&1&1&2 \\
&0&1&3
\end{array} 
\right] 
\]

\[
B_4=\left[ 
\begin{array}{c|ccc}
&1&0&1 \\
I_3&1&1&2 \\
&0&1&4
\end{array} 
\right] 
B_5=\left[ 
\begin{array}{c|ccc}
&1&0&1 \\
I_3&1&1&4 \\
&0&1&2
\end{array} 
\right] 
B_6=\left[ 
\begin{array}{c|ccc}
&1&0&1 \\
I_3&1&1&3 \\
&0&1&4
\end{array} 
\right] 
\] 
\normalsize

Matrices $B_1$ and $B_2$ are geometrically equivalent because $B_1$ can be transformed to $B_2$ as follows:   
multiply row 2 by 4; swap rows 1 and 3; replace row 2 by $row_2+row_1$; multiply column 2 by 4; swap columns 1 and 5 and columns 3 and 4. Similarly, we can check that matrix $B_1$ is geometrically equivalent to $B_4$ and matrix $B_2$ is geometrically equivalent to $B_3$ and $B_5$.
An exhaustive search found no mapping from $B_1$ to $B_3$ using the permissible operations. Thus we conclude that $Q_6$ has two geometrically inequivalent representations over $GF(5)$, $B_1$ and $B_3$. 

To confirm our exhaustive search we calculated the single-element extensions of $B_1$ and $B_3$.  The matroid $M[B_1]$ has 9 non-isomorphic simple single-element extensions over $GF(5)$ whereas the matroid $M[B_3]$ has 7 non-isomorphic simple single-element extensions. The first six matroids  in Figure 5 are common single-element extensions of $M[B_1]$ and $M[B_3]$. The seventh matroid is a single-element extension of $M[B_3]$, but not of $M[B_1]$. The last three matroids are single-element extensions of $M[B_1]$, but not of $M[B_3]$.

\begin{figure}[h]
\centering
\epsfxsize 6in \epsfbox{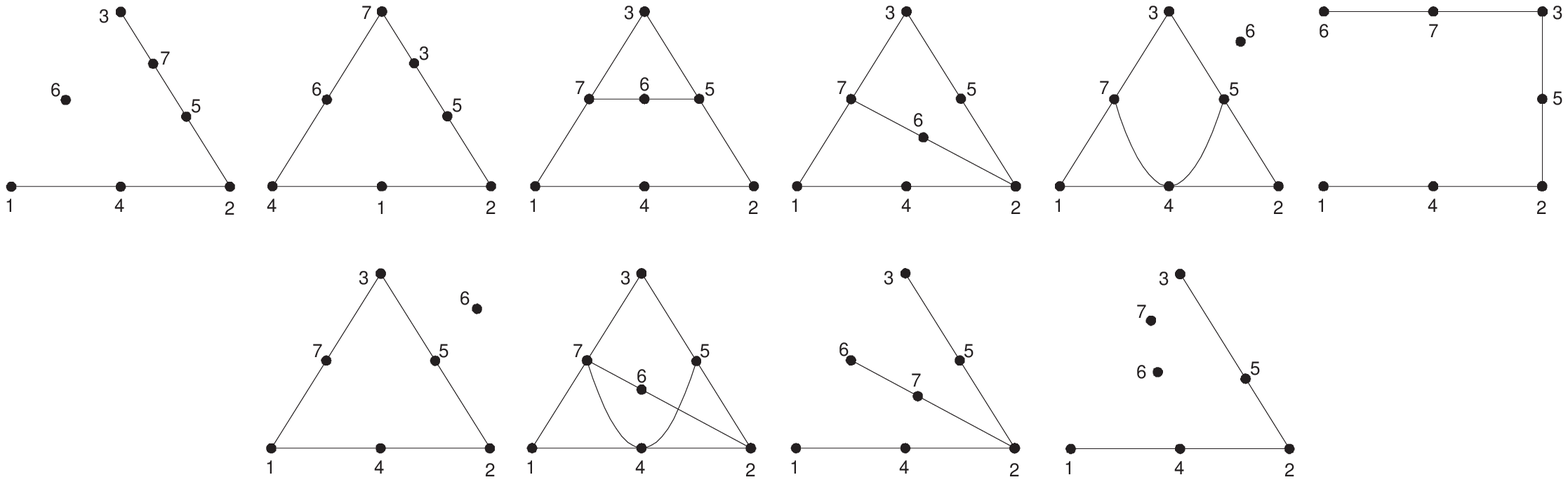}
\caption{Single-element extensions of $Q_6$ over $GF(5)$ \label{ExtnQ6}}
\end{figure}
\bigskip

As another example, we can check that $P_7$ has three projectively inequivalent representations over $GF(5)$, but is uniquely representable with respect to geometric equivalence. Whereas, over $GF(7)$  it has two geometrically inequivalent representations. The next two results follow from the above discussion.

\bigskip
\noindent {\bf Proposition 2.2.}  {\it If $A$ and $A'$ are projectively equivalent representations of $M$ over $F$, then they are geometrically equivalent.} $\qed$
\bigskip

As a result, if $M$ is uniquely representable over $GF(q)$ with respect to projective equivalence, then it is also uniquely representable over $GF(q)$ with respect to geometric equivalence. The converse of Proposition 2.2 is clearly not true since $W^3$ has three projectively inequivalent representations over $GF(5)$, but two of them are geometrically equivalent.
\bigskip

\noindent {\bf Proposition 2.3.}  {\it If $A$ and $A'$ are geometrically equivalent representations of a matroid $M$ over $F$, then $A$ and $A'$ have the same isomorphism classes of single-element extensions over $F$.}$\qed$
\bigskip

This raises the question ``Is there a one-one correspondence between existence of a linear transformation and having exactly the same isomorphism classes of single-element extensions?" The answer is no. For example, $W^3$ over $GF(7)$ has three geometrically inequivalent representations, but two of them have exactly the same single-element extensions.

\section {Geometric Equivalence versus stability}

This section gives a method for generating single-element extensions using geometric equivalence. A representation of $PG(r-1, q)$ in standard form may be written as $P=[I_r|L]$ 
by selecting from each one-dimensional subspace of $V(r, q)$, the column with 
one in the first non-zero position.  Let $M$ be a connected rank $r$ simple matroid with ground set $E$ and let 
$[I_r|D]$ be a representation of $M$ over $GF(q)$ in standard form.     
 If $[I_r|D]$ is a matrix over $GF(q)$ representing  $M$, then $[I_r|D|x]$ represents $M+e$, where $x$ is the column in the matrix $P$ corresponding to element $e$. Observe that there may be several such columns all representing the same element $e$.

\begin {enumerate}
\item [(i)] For each single-element extension $M+e$, let $X$ be a set of columns in $L-D$ such that for every $x\in X$, $[I_r|D|x]$ is an $GF(q)$-representation of $M+e$.  In other words $X$ is the set of columns that give the same single-element extension up to isomorphism. 
\item [(ii)] Group the columns in $X$ as follows:  $x_1$ and $x_2$ are in the same group if they have non-zero entries in the same positions. 
\item [(iii)] Let $A_1=[I_r|D|x_1]$ and $A_2=[I_r|D|x_2]$, where $x_1$ and $x_2$ are  columns from the same group. Check if  $M[A_1]= M[A_2]$ by listing the circuits and comparing them.  
\item [(iv)] Applying Proposition 2.1, $A_1$ and $A_2$ are projectively inequivalent representations if and only if $M[A_1]=M[A_2]$. However, if a linear transformation can be found between $A_1$ and $A_2$, then they are geometrically equivalent. Keep only one member in each equivalence class of geometrically equivalent representations. 
\item [(v)] Repeat for each geometrically inequivalent representation of $M$. It follows from Proposition 2.3 implies that this method will give all the non-isomorphic single-element extensions of $M$.  
\end {enumerate}

 The above algorithm is illustrated by computing the single-element extensions of $F_7^-$ over $GF(5)$ with respect to geometric equivalence.  This same example is presented in [2] as an example of how to compute projectively inequivalent extensions. Although not necessary for understanding this paper, the reader may find it interesting to compare and contrast [2] on projective equivalence with this paper on geometric equivalence, to see how the computational perspective naturally leads to the development of the concept of geometric equivalence.

Matrix representations for $F_7^-$ and $PG(2, 5)$ are given below:

\small
\[ 
A=\left[ 
\begin{array}{ccc|cccc}
&&&1&0&1&1 \\
&I_3&&1&1&0&1 \\
&&&0&1&1&1
\end{array} 
\right] 
\]
\normalsize

\tiny
\[
P=\left[ 
\begin{array}{ccc|cccccccccccccccccccccccccccc}
&&&    0&0&0&0&1&1&1&1&1&1&1&1&1&1&1&1&1&1&1&1&1&1&1&1&1&1&1&1 \\
&I_3&& 1&1&1&1&0&0&0&0&1&1&1&1&1&2&2&2&2&2&3&3&3&3&3&4&4&4&4&4 \\
&&&    1&2&3&4&1&2&3&4&0&1&2&3&4&0&1&2&3&4&0&1&2&3&4&0&1&2&3&4
\end{array} 
\right] 
\] 

\normalsize

\noindent Comparing $A$ and the matrix $P$ representing $PG(2, 5)$ and adding to $A$ the columns in $P$ missing in 
$A$ gives us three isomorphism classes:

\small
\begin{description}
\item[Class 1] $\{\{[0, 1, 2], [0, 1, 3]\}, \{[1, 2, 0], [1, 3, 0]\}, 
\{[1, 0, 2], [1, 0, 3]\}, $ \break 
$\{[1, 1, 3], [1, 1, 4], [1, 2, 2], [1, 3, 1], [1, 4, 1], [1, 4, 4]\}\}$
\item[Class 2] $\{\{[0, 1, 4]\}, \{[1, 4, 0]\}, \{[1, 0, 4]\}, \{[1, 1, 2], [1, 
2, 1], [1, 3, 3]\} \}$
\item[Class 3] $\{\{[1, 2, 3], [1, 2, 4], [1, 3, 2], [1, 3, 4], [1, 4, 2], [1, 
4, 3]\}\}$
\end{description}
\normalsize

\noindent Within each isomorphism class, group together the columns with non-zero entries in the 
same position and check if the resulting matroids are equal.  The three matroids corresponding to the 
three isomorphism classes, $M_1$, $M_2$, and $M_3$ are shown in Figure 4.  The new element is circled in each matroid.

\begin{figure}[h]
\centering
\epsfxsize 3.5in \epsfbox{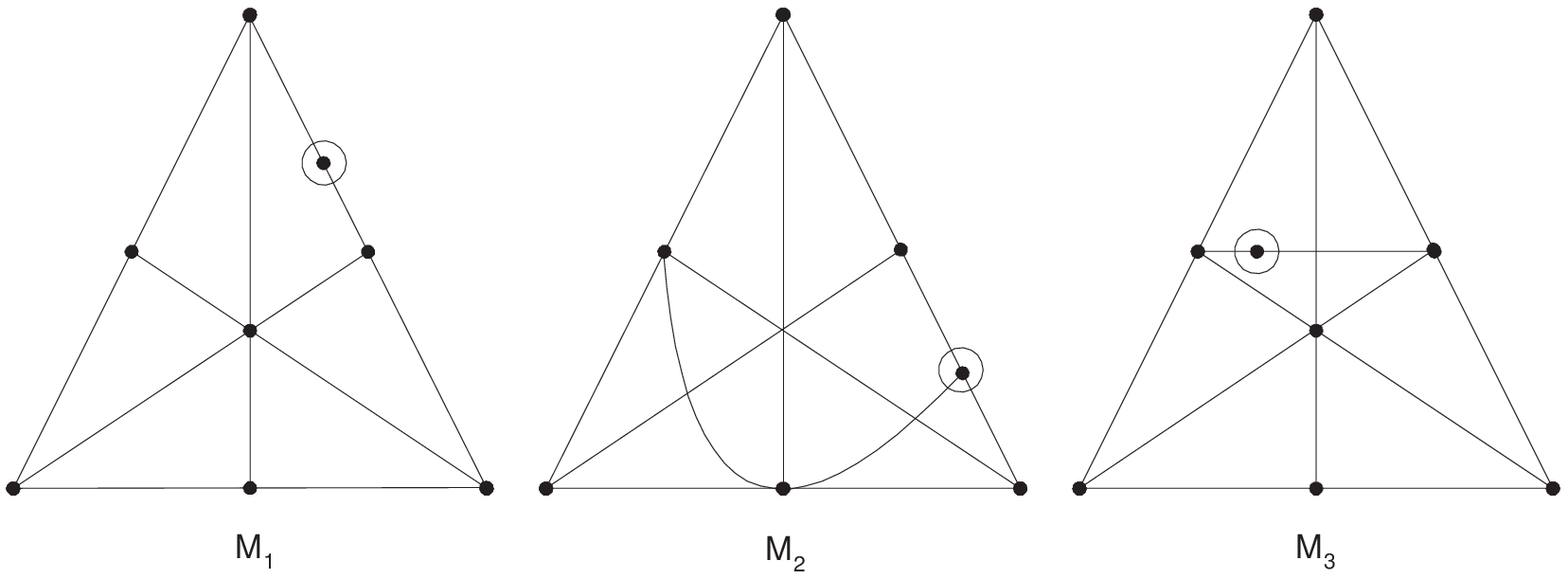}
\caption{The simple single-element extensions of $F_7^-$ over 
$GF(5)$ \label{SEE}}
\end{figure}
 
  In the first isomorphism class the matroid obtained by 
adding column $[0, 1, 2]$ to $A$ is equal to the matroid obtained by 
adding column $[0, 1, 3]$ to $A$.   However, $A$ with column $[0, 1, 2]$ can be transformed into $A$ with column $[0, 1, 3]$ as follows: swap rows 2 and 3; multiply rows 2 and 3 by 4; replace row 2 by $row_2+row_1$; replace row 3 by $row_3+row_1$; multiply columns 2, 3, and 5 by 4 and column 8 by 2; and swap columns 1 and 7. Thus these two representations are projectively inequivalent, but geometrically equivalent.
Similarly the matroid obtained by adding column $[1, 2, 0]$ to $A$ is equal to the matroid obtained by adding column $[1, 3, 0]$ to $A$, and the corresponding matrices are geometrically equivalent. 
The matroid obtained by adding column $[1, 0, 2]$ to $A$ is equal to the matroid obtained by adding column $[1, 0, 3]$ to $A$ and the two matrices are geometrically equivalent. This takes care of the first three groups in Class 1.
The situation is a bit more complicated in the last group consisiting of  Class 1 containing columns
$$\{[1, 1, 3], [1, 1, 4], [1, 2, 2], [1, 3, 1], [1, 4, 1], [1, 4, 4]\}\}$$
Here the matroid obtained by adding column $[1, 1, 3]$ to $A$ is equal to the matroid obtained by adding column $[1, 1, 4]$ to $A$.
The matroid obtained by adding column $[1, 2, 2]$ to $A$ is equal to the matroid obtained by adding column $[1, 4, 4]$ to $A$.
The matroid obtained by adding column $[1, 3, 1]$ to $A$ is equal to the matroid obtained by adding column $[1, 4, 1]$ to $A$. 
But the three matroids, $A$ with columns $[1, 1, 3]$, $[1, 2, 2]$, and $[1, 3, 1]$, respectively, are not equal (see Figure 5). We can check that in each of the three cases the corresponding matrices are geometrically equivalent.

Thus we may conclude that $F_7^-$ extends to  $M_1$ by giving rise to two projectively inequivalent representations of $M_1$. However, it extends to $M_1$ without giving rise to more geometrically inequivalent representations.

\begin{figure}[h]
\centering
\epsfxsize 3in \epsfbox{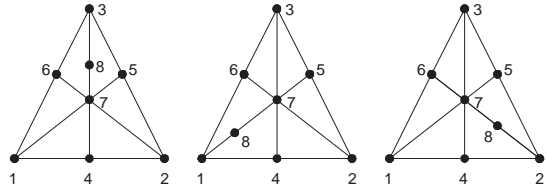}
\caption{Drawings of the first isomorphism class of $F_7^-$ over $GF(5)$ \label{ext2}}
\end{figure}

For the second isomorphism class, the three matroids, $A$ with columns $[1, 1, 2]$, $[1, 2, 1]$, and $[1, 3, 3]$, respectively, are not equal (see Figure 6).  So $F_7^-$ extends to $M_2$ without giving rise to any projectively inequivalent extensions of $M_2$ and consequently without giving rise to geometrically inequivalent representations of $M_2$.

\begin{figure}[h]
\centering
\epsfxsize 3.5in \epsfbox{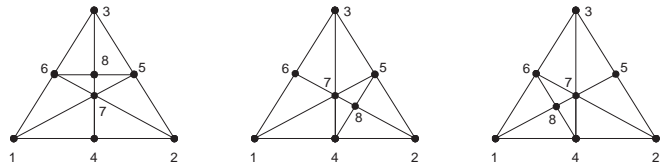}
\caption{Drawings of the second isomorphism class of $F_7^-$ over $GF(5)$  \label{ext2}}
\end{figure}

For the third isomorphism class, the matroid obtained by adding column $[1, 2, 3]$ to $A$ is equal to the matroid obtained by adding column $[1, 3, 4]$ to $A$, and the corresponding matrices are geometrically equivalent. The matroid obtained by adding column $[1, 2, 4]$ to $A$ is equal to the matroid obtained by adding column $[1, 4, 2]$ to $A$, and the corresponding matrices are geometrically equivalent. The matroid obtained by adding column $[1, 3, 2]$ to $A$ is equal to the matroid obtained by adding column $[1, 4, 3]$ to $A$, and the corresponding matrices are geometrically equivalent. Figure 7 shows that the three matroids, $A$ with columns $[1, 2, 3]$, $[1, 2, 4]$, and $[1, 3, 2]$, respectively, are not equal. So $F_7^-$ extends to $M_3$ by giving rise to two projectively inequivalent representations of $M_3$. However, it extends to  $M_3$ without giving rise to more geometrically inequivalent representations of $M_3$.  

\begin{figure}[h]
\centering
\epsfxsize 3.5in \epsfbox{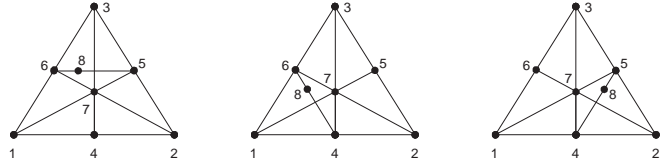}
\caption{Drawings of the third isomorphism class of $F_7^-$ over $GF(5)$  \label{ext3}}
\end{figure}

Moreover, since $F_7^-$ is uniquely representable over $GF(5)$ with respect to geometric equivalence, only one representation of $F_7^-$ may be used to get all its single-element extensions up to geometric equivalence.

\bigskip
\noindent {\bf Proposition 3.1.} {\it Let $M$ be a connected simple matroid represented by $A=[I_r|D]$   over $GF(q)$. Let $A_1=[I_r|D|x_1]$ and $A_2=[I_r|D|x_2]$ be matrices representing the same single-element extension $M+e$ where $x_1$ and $x_2$ are distinct columns in $P-D$ with non-zero entries in the same position. Then 

\begin{enumerate}
\item [(i)] $A_1$ and $A_2$ are projectively inequivalent representations of $M+e$ if and only if $M[A_1]=M[A_2]$.
\item [(ii)] $A_1$ and $A_2$ are geometrically inequivalent representations of $M+e$ if and only if $M[A_1]=M[A_2]$ and there is no linear transformation from $A_1$ to $A_2$. $\qed$
\end {enumerate}
} 
\bigskip

Note that, although the computations appear tedious to do by hand, from an algorithmic standpoint, geometric equivalence is easy in the sense that an algorithm for checking for linear transformations is a polynomial time algorithm.  The next proposition follows from the above discussion.

The concept of stability is a way of handling the large number of projectively inequivalent representations for a representable matroid. A minor $N$ of a matroid $M$ stabilizes $M$ over $GF(q)$ if no $GF(q)$-representation of $N$ can be extended to two projectively inequivalent $GF(q)$-representations of $M$. Observe that if $N$ has $k$ projectively inequivalent representations, then $M$ has at most $k$ projectively inequivalent representations. We say that $N$ is a {\it stabilizer} for $GF(q)$ if $N$ stabilizes each 3-connected $GF(q)$-representable matroid that contains $N$ as a minor. In the previous example, $F_7^-$ stabilizes its second extension, but does not stabilize its first and third extension.  
Contrast with the fact that $F_7^-$ extends to its single-element extensions without giving rise to new geometrically inequivalent representations. 

There is no analog of  the concept of stability for geometric equivalence.  For example, as mentioned earlier $W^3$ has two inequivalent geometric representations over $GF(5)$. One of its single-element $B_7$ is uniquely representable with respect to geometric equivalence. However, a single-element extension of $B_7$ has three geometric inequivalent representations over $GF(5)$.  Geometric equivalence may be viewed as an  alternative to stability. Future research includes finding bounds on the number of elements in a rank $r$ matroid beyond which all the matroids are uniquely representable with respect to geometric equivalence.

\bigskip

\noindent {\bf References}

\begin{enumerate}

\item T. H. Brylawski and D. Lucas (1976). Uniquely representable combinatorial geometries. In {\it Teorie combinatorie (Proc. 1973 Internat. Colloq.)}, (1976), pp. 83-104. Accademia Nazionale dei Lincei, Rome. 

\item J. Geelen, G. Whittle (2010). The projective plane is a stabilizer, {\it J. of Combin. Theory Ser. B}, {\bf 100}, no. 2  128-131. 

\item S. R. Kingan (2009). A computational approach to inequivalence in matroids. Proceedings of the Thirty-Seventh Southeastern International Conference on Combinatorics, Graph Theory and Computing. {\it Congr. Numer.} {\bf 198}, 63-74.

\item J. G. Oxley (1992). {\it Matroid theory}, (1992), Oxford University Press, New York. 

\end{enumerate}

\end {document}